\documentclass[leqno,draft]{article}
\usepackage{bcp04e}
\usepackage{amsmath}
\usepackage{amsfonts}

\newcommand{\field}[1]{\mathbb{#1}}

\newcommand{\RR}{\field{R}}

\def\Ad{\mathop{\rm Ad}\nolimits}

\def\div{\mathop{\rm div}\nolimits}
\def\rot{\mathop{\rm rot}\nolimits}
\def\GL{\mathop{\rm Gl}\nolimits}
\def\Aff{\mathop{\rm Aff}\nolimits}
\def\P{\mathop{\rm P}\nolimits}
\def\PGL{\mathop{\rm PGl}\nolimits}
\def\O{\mathop{\rm O}\nolimits}
\def\U{\mathop{\rm U}\nolimits}

\def\SU{\mathop{\rm SU}\nolimits}
\newdefin{dfn}{Definition}
\newtheorem{theorem}{Theorem}
\newtheorem{proposition}[theorem]{Proposition}
\newdefin{remark}{Remark}

\begin{document}
\mathclass{Primary53C05; Secondary 53B05, 53B10, 53B15, 70G45.}

\abbrevauthors{C.-M.~Marle} \abbrevtitle{The works of
Charles~Ehresmann on connections}

\title{The works of Charles Ehresmann on connections:\\
from Cartan connections\\ to connections on fibre bundles}

\author{Charles-Michel Marle}
\address{Institut de Math\'{e}matiques, Universit\'{e} Pierre et Marie Curie,\\
4, place Jussieu, 75252 Paris c\'{e}dex 05, France\\
E-mail: charles-michel.marle@polytechnique.org}
\thanks{}

\maketitlebcp

\begin{abstract}
Around 1923, \'Elie Cartan introduced affine connections on
manifolds and defined the main related concepts: torsion, curvature,
holonomy groups. He discussed applications of these concepts in
Classical and Relativistic Mechanics; in particular he explained how
parallel transport with respect to a connection can be related to
the principle of inertia in Galilean Mechanics and, more generally,
can be used to model the motion of a particle in a gravitational
field. In subsequent papers, \'Elie Cartan extended these concepts
for other types of connections on a manifold: Euclidean, Galilean
and Minkowskian connections which can be considered as special types
of affine connections, the group of affine transformations of the
affine tangent space being replaced by a suitable subgroup; and more
generally, conformal and projective connections, associated to a
group which is no more a subgroup of the affine group.
\par\smallskip
Around 1950, Charles Ehresmann introduced connections on a fibre
bundle and, when the bundle has a Lie group as structure group,
connection forms on the associated principal bundle, with values in
the Lie algebra of the structure group. He called {\it Cartan
connections\/} the various types of connections on a manifold
previously introduced by \'E.~Cartan, and explained how they can be
considered as special cases of connections on a fibre bundle with a
Lie group $G$ as structure group: the standard fibre of the bundle
is then an homogeneous space $G/G'$; its dimension is equal to that
of the base manifold; a Cartan connection determines an isomorphism
of the vector bundle tangent to the the base manifold onto the
vector bundle of vertical vectors tangent to the fibres of the
bundle along a global section.
\par\smallskip
These works are reviewed and some applications of the theory of
connections are sketched.
\end{abstract}

\section{Introduction.}

It is a great honor for me to be invited at the Seventh
International Conference on Geometry and Topology of Manifolds,
dedicated to the mathematical legacy of Charles Ehresmann. I enjoyed
with great pleasure the hospitality of the Mathematical Research and
Conference Center of the Polish Academy of Sciences, and I address
my warmest thanks to the organizers and to the supporting
institutions.

\par\smallskip
Around 1923, \'Elie Cartan~\cite{cartan1, cartan2, cartan3}
introduced the notion of an \emph{affine connection} on a manifold.
That notion was previously used, in a less general setting, by
H.~Weyl~\cite{weyl1} and rests on the idea of parallel transport due
to T.~Levi-Civita~\cite{levicivita}.
\par\smallskip
At the very beginning of \cite{cartan1}, even before defining
explicitly affine connections, \'Elie Cartan explains how that
concept can be used in Newtonian and Einsteinian Mechanics. He shows
that the \emph{principle of inertia} (which is at the foundations of
Mechanics), according to which a material point particle, when no
forces act on it, moves along a straight line with a constant
velocity, can be expressed locally by the use of an affine
connection. Under that form, that principle remains valid in
(curved) Einsteinian  space-times.
\par\smallskip
Cartan even shows that by a suitable adjustment of the connection, a
gravity force (that means, an acceleration field) can be accounted
for, and becomes a part of the Geometry of space-time. That result
expresses the famous \emph{equivalence principle} used by Einstein
for the foundations of General Relativity. As shown by Cartan, it is
valid for Newtonian Mechanics as well.
\par\smallskip
Then \'Elie Cartan presents a thorough geometric study of affine
connections; he defines their curvature and torsion, and discusses
the parallel displacement of a frame along a closed loop. He
introduces Euclidean, Galilean and Minkowskian connections, for
which the group of affine transformations is replaced by a suitable
subgroup. In \cite{cartan4, cartan5, cartan6} he introduces more
general types of connections associated to transformation groups
which are no more subgroups of the group of affine transformations.
\par\smallskip
Cartan's ideas were fully formalized by Charles Ehresmann in the
framework of connections on fibre bundles, which he introduced in
\cite{ehresmann2}.
\par\smallskip
In Section~2 we briefly present Cartan's intuitive ideas about
connections. Then in Section~3 we describe Ehresmann connections on
fibre bundles, and in Section~4 Cartan connections as seen by
Ehresmann. In Section~5 we present with more details examples of
Cartan connections,  including affine, projective and conformal
connections. In Section~6, following Cartan, we show how a
gravitational force can be included in the geometry of space-time by
the use of a suitable connection, and we briefly present other
applications of connections: Geometric quantization, phases in
mechanics, nonholonomic or active constraints, Maxwell's equations,
Yang-Mills fields.

\section{Cartan affine connections and their generalizations.}

\'E.~Cartan writes in \cite{cartan1}: \lq\lq Une vari\'{e}t\'{e} \`{a} connexion
affine est une vari\'{e}t\'{e} qui, au voisinage imm\'{e}diat de chaque point, a
tous les caract\`{e}res d'un espace affine, et pour laquelle on a une
loi de rep\'{e}rage des domaines entourant deux points infiniment
voisins~: cela veut dire que si, en chaque point, on se donne un
syst\`{e}me de coordonn\'{e}es cart\'{e}siennes ayant ce point pour origine, on
conna\^{\i}t les formules de transformation (de m\^{e}me nature que dans
l'espace affine) qui permettent de passer d'un syst\`{e}me de r\'{e}f\'{e}rence
\`{a} tout autre syst\`{e}me de r\'{e}f\'{e}rence d'origine infiniment
voisine\rq\rq.
\par\smallskip
Approximate translation:
\lq\lq A manifold with an affine connection
is a manifold whose properties, in the neighborhood of each point,
are those of an affine space, and on which there is a law for
fitting together the neighborhoods of two infinitesimally nearby
points: it means that if, on a neighborhood of each point, we have
chosen Cartesian coordinates with that point as origin, we know the
transformation formulae (of the same nature as those valid in an
affine space) which allow to go from a reference frame to another
reference frame with an infinitesimally nearby origin\rq\rq.
\par\smallskip
In fact, given a smooth  manifold $M$, it is not on a neighborhood
of each point $m\in M$ that Cartan considers a local affine
structure. Rather, at each point $m\in M$, he considers the tangent
space $T_mM$ endowed with its natural affine space structure. And he
writes:
\lq\lq La vari\'{e}t\'{e} sera dite \`{a} \emph{connexion affine}
lorsqu'on aura d\'{e}fini, d'une mani\`{e}re d'ailleurs arbitraire, une loi
permettant de rep\'{e}rer l'un par rapport \`{a} l'autre les espaces affines
attach\'{e}s \`{a} deux points infiniment voisins quelconques $m$ et $m'$ de
la vari\'{e}t\'{e}; cette loi permettra de dire que tel point de l'espace
affine attach\'{e} au point $m'$ correspond \`{a} tel point de l'espace
affine attach\'{e} au point $m$, que tel vecteur du premier espace est
parall\`{e}le ou \'{e}quipollent \`{a} tel vecteur du second espace. En
particulier le point $m'$ lui-m\^{e}me sera rep\'{e}r\'{e} par rapport \`{a}
l'espace affine du point $m$ $\ldots$\rq\rq.
\par\smallskip

Approximate translation:
\lq\lq The manifold will be said to be
endowed with an \emph{affine connection} once we have defined, in an
arbitrary way, a law allowing to localize one with respect to the
other the affine spaces attached to two infinitesimally nearby
points $m$ and $m'$ of that manifold; that law will tell us which
point of the affine space attached to $m'$ corresponds to a given
point of the space attached to $m$, and will tell us whether a
vector living in the first space is parallel, or equipollent, to a
vector living in the second space; in particular, the point $m'$
itself will be localized in the affine space attached to $m$
$\ldots$\rq\rq.
\par\smallskip
Cartan then explains with more details how to define an affine
connection on a $3$-dimensional manifold $M$; of course his
definition extends easily to manifolds of any dimension. He
considers, at each point $m\in M$, an affine frame of the affine
tangent space $T_mM$, with as origin the point $m$ itself
(identified with the null vector at $m$), and with the linear basis
$(e_1,e_2,e_3)$ as basis. In order to define the law which links the
affine spaces tangent to the manifold $M$ at two infinitesimally
nearby points $m$ and $m'$, he write the relations
\begin{equation}
\begin{split}
 dm&=\omega^1e_1+\omega^2e_2+\omega^3e_3\,,\\
 de_i&=\omega^1_ie_1+\omega^2_ie_2+\omega^3_ie_3\,,\quad i=1,\ 2,\ 3.
 \end{split}
 \end{equation}
\par\smallskip
These equations mean that the point $m'$, origin of $T_{m'}M$,
infinitesimally near $m$, must be identified with the point
 $$m+\omega^1e_1+\omega^2e_2+\omega^3e_3$$
of the affine space $T_mM$.
\par\smallskip
Similarly, the vectors $e'_1$, $e'_2$, $e'_3$ of $T_{m'}M$ must be
identified with the vectors
  $$e'_i=e_i+\omega^1_ie_1+\omega^2_ie_2+\omega^3_ie_3\,,\ i=1,\ 2,\ 3$$
of $T_mM$.

\par\smallskip
Equations $(1)$ must be understood as equalities between
differential $1$-forms on a $15$-dimensional space with, as
coordinates, the $3$ coordinates which specify a point on $M$, and
$12$ more coordinates on which depend the affine frames of a
$3$-dimensional affine space. In fact, these differential $1$-forms
live on the principal bundle of affine frames of the affine tangent
spaces to the maniflold $M$. The meaning of Equations $(1)$ is the
following: let $h=(m,e_1,e_2,e_3)$ be an affine frame of the affine
tangent space $T_mM$ and let $V$ be a vector tangent at $h$ to the
space of affine frames. The infinitesimal displacement of that
affine frame defined by the vector $V$ is parallel, with respect to
the connection, if and only if $\omega(V)=0$, where $\omega$ is the
$1$-form of the connection (in the sense of Ehresmann, as introduced
later). Equation $\omega(V)=0$, when explicitly written, leads to
the $12$ scalar equations $(1)$ because the $1$-form $\omega$ takes
its values in the $12$-dimensional Lie algebra of the group of
affine transformations of ${\RR}^3$. This will become clear with the
works of Charles Ehresmann~\cite{ehresmann2}.
\par\smallskip
Euclidean connections, Galilean connections, Minkowskian
connections, $\ldots$, considered by Cartan, appear as special cases
of affine connections, obtained by replacing the group of affine
transformations of an affine space by an appropriate subgroup. For
example, on a $3$-dimensional manifold $M$ endowed with a Riemannian
metric $g$, instead of general affine frames $h=(m,e_1,e_2,e_3)$
with the contact point $m$ as origin of the affine tangent space
$T_mM$, Cartan considers only orthonormal frames, which satisfy
 $$g(e_i,e_j)=\delta_{i\,j}\,.$$
Now the dimension of the space of orthonormal affine frames is 9
(instead of 15 for the space of general affine frames). By
differentiation of the above relations we see that the $1$ forms
$\omega^i$ and $\omega^j_i$ of Equations $(1)$ must satisfy
 $$\omega^1_1=\omega^2_2=\omega^3_3=0\,,\quad \omega^3_2+\omega^2_3=0\,,
 \quad \omega^1_3+\omega^3_1=0\,,
 \quad \omega^2_1+\omega^1_2=0\,.
 $$
Only 6 of the 12 scalar equations $(1)$ are now independent (in
agreement with the dimension, 6,  of the Lie algebra of the group of
affine orthogonal transformations of ${\RR}^3$).
\par\smallskip
Other types of connections were defined and discussed by
\'E.~Cartan~\cite{cartan4, cartan5, cartan6}, with transformation
groups which are no more subgroups of the group of affine
transformations: conformal connections, projective connections,
$\ldots$ The general idea underlying the notions of such connections
on a manifold $M$ is to use, as a local model of that manifold, an
homogeneous space of the same dimension as $M$. A copy of that
homogeneous space is attached at each point of the manifold, and
considered as \lq\lq tangent\rq\rq\ at that point to the manifold.
The connection is essentially a law which indicates how these
homogeneous spaces are glued together.
\par\smallskip
Of course, affine connections appear as the special case in which
the homogeneous space which is attached to each point of the
manifold is an affine space, more precisely the affine space tangent
to the manifold at that point.
\par\smallskip
For conformal connections on an $n$-dimensional manifold,
\'E.~Cartan~\cite{cartan4} writes: \lq\lq Attachons \`{a} chaque point
$P$ de cette vari\'{e}t\'{e} un espace
 conforme \`{a} $n$ dimensions, $\ldots$ La vari\'{e}t\'{e} sera dite \emph{\`{a}
 connexion conforme} si nous nous donnons une loi (d'ailleurs
 arbitraire) permettant de rapporter, d'une mani\`{e}re conforme,
 l'espace conforme attach\'{e} au point $P$ de la vari\'{e}t\'{e} \`{a} l'espace
 conforme attach\'{e} au point infiniment voisin $P'\,\,$\rq\rq.
 \par\smallskip

Approximate translation: \lq\lq Let us link an $n$-dimensional
conformal space to each point $P$ of our manifold. That manifold
will be said to be endowed with a \emph{conformal connection} when
we have specified, in an arbitrary way, how to tie (or maybe
identify) the conformal space linked at point $P$ with the conformal
space linked to the infinitesimally nearby point $P'\,\,$\rq\rq.
\par\smallskip

In~\cite{cartan4}, Cartan writes:\quad
 \lq\lq L'id\'{e}e fondamentale se rattache \`{a} la notion de
 \emph{parall\'{e}lisme} que M.~T.~Levi-Civita a introduite de mani\`{e}re
 si f\'{e}conde. Les nombreux auteurs qui ont g\'{e}n\'{e}ralis\'{e} la th\'{e}orie
 des espaces m\'{e}triques sont tous partis de l'id\'{e}e fondamentale de
 M.~Levi-Civita, mais, semble-t-il, sans pouvoir la d\'{e}tacher de
 l'id\'{e}e de \emph{vecteur}. Cela n'a aucun inconv\'{e}nient quand il
 s'agit de vari\'{e}t\'{e}s \`{a} connexion affine $\ldots$ Mais cela semblait
 interdire tout espoir de fonder une th\'{e}orie \emph{autonome} de
 vari\'{e}t\'{e}s \`{a} connexion conforme ou projective. En fait, ce qu'il y
 a d'essentiel dans l'id\'{e}e de M.~Levi-Civita, c'est qu'elle donne
 un moyen pour raccorder entre eux deux petits morceaux infiniment
 voisins d'une vari\'{e}t\'{e}, et c'est cette id\'{e}e de \emph{raccord}
qui est f\'{e}conde\rq\rq.
\par\smallskip
Approximate translation: \lq\lq The fundamental idea stems from the
notion of \emph{parallelism} introduced by M.~T.~Levi-Civita in such
a fruitful way. The many authors who genralized the theory of metric
spaces all started from the fundamental idea of M.~Levi-Civita, but,
seemingly, without freeing it from the idea of \emph{vector}. That
does not matter as long as one deals with manifolds with affine
connections $\ldots$ But that seemed to forbid any hope to build an
\emph{autonomous} theory of manifolds with conformal or projective
connections. In fact, the main thing in M.~Levi-Civita's idea is
that it allows to glue together two small, infinitesimally nearby
pieces of a manifold, and it is that idea of gluing which is most
fruitful\rq\rq.

\section{Ehresmann connections.}

Let $E(B,F)$ be a locally trivial smooth fibre bundle with base $B$,
standard fibre $F$, and canonical projection $\pi_E:E\to B$. For
each $x\in B$, the fibre at $x$, $E_x=\pi_E^{-1}(x)$, is
diffeomorphic to $F$. Ehresmann~\cite{ehresmann2} defines an
\emph{infinitiesimal connection} on that bundle as follows.

\begin{dfn}\label{definfinitesimalconnection} An \emph{infinitesimal connection}
on the locally trivial bundle $E(B,F)$ is a vector sub-bundle $C$ of
$TE$, complementary to $\ker(T\pi_E)$, {\it i.e.} such that for each
$z\in E$,
 $$T_zE=\ker(T_z\pi_E)\oplus C_z\,,$$
which satisfies the additional condition:
 \par\smallskip
$(C)$\quad Given any smooth path $t\mapsto x(t)$ in $B$ going from a
point $x_0=x(t_0)$ to another point $x_1=x(t_1)$, and any $z_0\in
E_{x_0}$, there exists a smooth path $t\mapsto z(t)$ in $E$, called
the \emph{horizontal lift} of $t\mapsto x(t)$ through $z_0$, such
that
 $$z(t_0)=z_0\,,\quad\pi_E\bigl(z(t)\bigr)=x(t)\ \hbox{and}\
 \frac{dz(t)}{dt}\in C_{z(t)}
 \quad\hbox{for all}\quad t\in[t_0,t_1]\,.$$
More generally any smooth path $t\mapsto z(t)$ in $E$ is said to be
\emph{horizontal} if for each $t$, $\displaystyle
\frac{dz(t)}{dt}\in C_{z(t)}$.
\end{dfn}

\begin{remark}\label{remparalleltransport}
Under the assumptions of Definition 1, let $t\mapsto x(t)$ be a
smooth path in $B$, defined on an interval $I$. For each $t_0\in I$
and each $z_0\in E_{x(t_0)}$, its horizontal lift $t\mapsto z(t)$
through $z_0$ is unique and defined on $I$. For each pair
$(t_0,t_1)\in I^2$, the map $\varphi_{t_1\,t_0}:E_{x(t_0)}\to
E_{x(t_1)}$,
 $$z_0\mapsto\varphi_{t_1\,t_0}(z_0)=z(t_1)\,,$$
is a diffeomorphism of the fibre $E_{x(t_0)}$ onto the fibre
$E_{x(t_1)}$, called the \emph{parallel transport} of the fibres of
$E(B,F)$ along the smooth path $t\mapsto x(t)$. By looking only at
paths defined on finite closed intervals, and taking their end
points, we see that the connection $C$ determines a homomorphism of
the groupoid of smooth paths in $B$ which joint two points in $B$
into the groupoid of diffeomorphisms of a fibre or $E(B,F)$ onto
another fibre.
\end{remark}

\begin{dfn}\label{defdevelopment}
Let $t\mapsto z(t)$ be a smooth path, defined on an interval $I$,
with values in the total space $E$ of a locally trivial, smooth
bundle $E(B,F)$ equipped with a connection $C$. Let $t\mapsto
x(t)=\pi_E\bigl(z(t)\bigr)$ be the corresponding smooth path in the
base manifold $B$. Let $t_0\in I$ and, for each $t$, let
$\varphi_{t\,t_0}$ be the parallel transport of the fibre
$E_{x(t_0)}$ onto the fibre $E_{x(t)}$ along the path $t\mapsto
x(t)$. The smooth path in the fibre $E_{x(t_0)}$,
  $$t\mapsto \varphi_{t\,t_0}^{-1}\circ z(t)$$
is called the \emph{development} of the path $t\mapsto z(t)$ in the
fibre $E_{x(t_0)}$.
\end{dfn}

\begin{remark} The development in $E_{x(t_0)}$ of a path $t\mapsto
z(t)$ in $E$ is a constant path if and only if $t\mapsto z(t)$ is
horizontal.
\end{remark}

We assume now that $E(B,F)$ is a smooth fibre bundle $E(B,F,G,H)$
with a Lie group $G$ as structure group. Let us first recall that
notion, due to Ehresmann, and explain the notations. For simplicity
we consider only smooth bundles, although in \cite{ehresmann1}
Ehresmannn, more generally, defines topological bundles. The
\emph{total space} $E$ and the \emph{base} $B$ of the bundle are two
smooth manifolds, and there is a smooth submersion $\pi_E:E\to B$ of
the total space onto the base. For each $x\in B$,
$E_x=\pi_E^{-1}(x)$ is a smooth submanifold of $E$ called the
\emph{fibre over} $x$. The \emph{standard fibre} $F$ is a smooth
manifold and the \emph{structure group} $G$ is a Lie group of
diffeomorphisms of $F$, {\it i.e.} a Lie group $G$ which acts on the
manifold $F$ on the left in such a way that the map $G\times F\to
F$, $(g,f)\mapsto gf$, is smooth. That action is assumed to be
effective, {\it i.e.} such that for each $g\in G$ other than the
unit element $e$, the corresponding diffeomorphism $f\mapsto gf$ of
the standard fibre $F$ is not the identity; therefore an element
$g\in G$ is fully determined by the corresponding diffeomorphism of
$F$. An \emph{atlas} $(U_i,\varphi_i\,,\ i\in I)$ of the bundle
$E(B,F,G,H)$ is defined by the following data:

\begin{enum}{(iii)}
\renewcommand\theenumi{\roman{enumi}}
\renewcommand\labelenumi{(\theenumi)}

\item a family $(U_i\,,\ i\in I)$ of open subsets of $B$ such that
$\bigcup_{i\in I}U_i=B$,

\item for each $i\in I$, a diffeomorphism
  $$\varphi_i:\pi_E^{-1}(U_i)\to U_i\times F\,,$$
such that for each $z\in\pi_E^{-1}(U_i)$,
  $$p_1\bigl(\varphi_i(z)\bigr)= \pi_E(z)$$
(where $p_1:U_i\times F\to U_i$ is the first projection) and that
for each pair $(i,j)\in I^2$ with $U_i\cap U_j\neq\emptyset$,
  $$\varphi_j\circ\varphi_i^{-1}(x\,,\,\xi)=\bigl(x\,,\,s_{j\,i}(x)(\xi)\bigr)\,,$$
with $x\in U_i\cap U_j$, $\xi\in F$, and where $s_{j\,i}:U_i\cap
U_j\to G$ is a smooth map which takes its values in the Lie group
$G$.
\end{enum}

Two atlases $(U_i,\varphi_i\,,\ i\in I)$ and $(V_j,\psi_j\,,\ j\in
J)$ of the bundle $E(B,F,G,H)$ are said to be \emph{equivalent} if
their union is still an atlas of that bundle. The $G$-bundle
structure of $E(B,F,G,H)$ is determined by the choice of an
equivalence class of atlases. Diffeomorphisms
$\varphi_i:\pi_E^{-1}(U_i)\to U_i\times F$ which belong to an atlas
of that equivalence class are called \emph{admissible local
trivializations} of the bundle $E(B,F,G,H)$.
\par\smallskip
For each $x\in B$, the fibre $E_x=\pi_E^{-1}(x)$ is diffeomorphic to
$F$. More precisely there are \emph{admissible diffeomorphisms} of
$F$ onto $E_x$, obtained by means of an admissible local
trivialization $(U_i,\varphi_i)$, with $x\in U_i$:
  $$\xi\mapsto \varphi_i^{-1}(x,\xi)\,,\quad\hbox{with}\ \xi\in
  F\,.$$
For each  $x\in E$, if $h:F\to\pi_E^{-1}(x)$ is an admissible
diffeomorphism, another diffeomorphism $h':F\to \pi_E^{-1}(x)$ is
admissible if and only if there exists $g\in G$ such that $h'=h\circ
g$.
\par\smallskip
In Ehresmann's notation $E(B,F,G,H)$, we have already explained what
are the total space $E$, the base $B$, the standard fibre $F$ and
the structure group $G$. We still have to define $H$: it is the set
of all admissible diffeomorphisms of the standard fibre $F$ onto the
fibre $E_x=\pi_E^{-1}(x)$, for all $x\in B$. In \cite{ehresmann1},
Ehresmann proves that $H$ is a smooth manifold. More precisely it is
the total space of a smooth fibre bundle $H(B,G,G_\gamma,\overline
H)$ with base $B$, standard fibre $G$, structure group $G_\gamma$
(that means the Lie group $G$ considered as acting on itself by left
translations). The fibre $H_x=\pi_H^{-1}(x)$ over a point $x\in B$
is the set of admissible diffeomorphisms of $F$ onto
$E_x=\pi_E^{-1}(x)$. In Ehresmann's notations $\overline H$ is the
set of admissible diffeomorphisms $\overline h$ of the standard
fibre $G$ onto the fibre $H_x$ of $H$ over some point $x\in B$. In
fact, $\overline H$ can be identified with $H$, since an element
$h\in H_x$, previously considered as a diffeomorphism of $F$ onto
$E_x$, can also be considered as a diffeomorphism $\overline h$ of
$G$ onto $H_x$ by writing
  $$\overline h(g)=h\circ g:F\to E_x\,.$$
The smooth fibre bundle $H(B,G,G_\gamma,\overline H)$ is said to be
a \emph{principal bundle} because its structure group $G$ is also
its standard fibre. More precisely, it is called \emph{the principal
bundle associated to} the bundle $E(B,F, G,H)$, and it is fully
determined by that bundle.
\par\smallskip
Conversely, let $H(B,G,G_\gamma,\overline H)$ be a principal bundle
with structure Lie group $G_\gamma$, and $F$ be a smooth manifold on
which the Lie group $G$ acts on the left. There exist a bundle
$E(B,F,G,H)$, unique up to an isomorphism, with structure Lie group
$G$, standard fibre $F$ and $H(B,G,G_\gamma,\overline H)$ as
associated principal bundle. The bundles obtained in that way are
said to be \emph{associated} to the given principal bundle
$H(B,G,G_\gamma,\overline H)$.

\begin{remark}\label{remrightaction}
There is a natural, free action on the right of the structure Lie
group $G$ on the total space $H$ of the principal bundle
$H(B,G,G_\gamma,\overline H)$,
  $$H\times G\to H\,,\quad (h,g)\mapsto R_g(h)=h\circ g\,,$$
whose orbits are the fibres $H_x$, $x\in B$.
\end{remark}

\begin{dfn}\label{defgconnection} Let $E(B,F,G,H)$
be a fibre bundle with a Lie group $G$ as structure group. An
infinitesimal connection $C$ on that bundle, in the sense of
Definition~\ref{definfinitesimalconnection}, is said to be
\emph{compatible with the structure group} $G$, and is called a
\emph{$G$-connection}, if the parallel transport along any smooth
path in $B$ starting from any point $x_0\in  B$ and ending at any
other point $x_1\in B$, is of the form $h_1\circ h_0^{-1}$, with
$h_0\in H_{x_0}$ and $h_1\in H_{x_1}$.
\end{dfn}

\begin{proposition}
A $G$-connection on the bundle $E(B,F,G,H)$ determines a unique
$G_\gamma$-connection $\overline C$ on the associated principal
bundle $H(B,G,G_\gamma,\overline H)$ such that, for any smooth path
$t\mapsto x(t)$ in the base $B$, any $t_0$ in the interval of
definition of that path and any $h_0\in H_{x(t_0)}$, the horizontal
lift through $h_0$ of the path $t\mapsto x(t)$ with respect to the
connection $\overline C$ is the path
 $$t\mapsto h(t)=\varphi_{t\,t_0}\circ h_0\,,$$
where $\varphi_{t\,t_0}:E_{x(t_0)}\to E_{x(t)}$ is the parallel
transport along the smooth path $t\mapsto x(t)$, for the
$G$-connection $C$ on $E(B,F,G,H)$.

Conversely, a $G_\gamma$-connection $\overline C$ on the principal
bundle $H(B,G,G_\gamma,\overline H)$ determines a unique
$G$-connection $C$ on $E(B,F,G,H)$ such that, for any smooth path
$t\mapsto x(t)$ in the base $B$, any $t_0$ in the interval of
definition of that path and any $z_0\in E_{x(t_0)}$, the horizontal
lift through $z_0$ of the path $t\mapsto x(t)$ with respect to the
connection $C$ is the path
 $$t\mapsto z(t)=h(t)(\xi)\,,$$
where $t\mapsto h(t)$ is any horizontal lift in $H$ of the path
$t\mapsto x(t)$ for the $G_\gamma$-connection $\overline C$, and
where $\xi=\bigl(h(t_0)\bigr)^{-1}(z_0)$ is the unique element in
$F$ such that $h(t_0)(\xi)=z_0$.
\end{proposition}

\begin{remark}
A $G_\gamma$-connection on a principal bundle
$H(B,G,G_\gamma,\overline H)$ is a vector sub-bundle $\overline C$
of the tangent bundle $TH$, complementary to $\ker T\pi_H$, {\it
i.e.} such that, for each $h\in H$,
  $$T_hH=\ker T_h\pi_H\oplus \overline C_h\,,$$
which remains invariant by the right action of $G_\gamma$ defined in
Remark~\ref{remrightaction}, {\it i.e.} which is such that for any
$h\in H$ and $g\in G$,
 $$TR_g(\overline C_h)=\overline C_{h\circ g}\,,$$
where $TR_g:TH\to TH$ is the lift to vectors of the map $R_g:H\to
H$, $h\mapsto h\circ g$. Conversely, a vector sub-bundle $\overline
C$ of $TH$ complementary to $\ker T\pi_H$ which remains invariant by
the right action of $G_\gamma$ is a $G$-connection on
$H(B,G,G_\gamma,\overline H)$: by using the $G$-invariance, one can
prove that condition $(C)$ of
Definition~\ref{definfinitesimalconnection} is automatically
satisfied.
\end{remark}

The next proposition introduces the very important notion of
\emph{connection form}.

\begin{proposition}\label{propconnectionform}
Let $\overline C$ be a $G_\gamma$-connection on the principal bundle
$H(B,G,G_\gamma,\overline H)$. There exists a unique $1$-form
$\omega$ on $H$, with values in the Lie algebra $\cal G$ of $G$,
such that for all $h\in H$,
 $$\overline C_h=\bigl\{\,V\in T_hH; \omega(V)=0\,\bigr\}\,,$$
and that for each $\eta\in \cal G$,
  $$\omega\bigl(\widehat \eta(h)\bigr)=\eta\,.$$
We have denoted by $\widehat \eta$ the fundamental vector field on
$H$ associated to $\eta$, {\it i.e.} the vector field whose value,
at an element $h\in H$, is
 $$\widehat\eta(h)=\frac{d\bigl(R_{\exp(t\eta)}h\bigr)}{dt}\Bigm|_{t=0}\,.$$
The $1$-form $\omega$ is called the connection form of the
$G_\gamma$-connection $\overline C$.
\end{proposition}

A connection on a principal bundle compatible with its structure
group can be defined by its connection form, as shown by the next
proposition.

\begin{proposition}\label{proppropertiesconnectionform}
The connection form $\omega$ of a $G_\gamma$-connection $\overline
C$ on the principal bundle $H(B,G,G_\gamma,\overline H)$ has the
following properties:

\begin{enum}{(iii)}

\renewcommand\theenumi{\roman{enumi}}
\renewcommand\labelenumi{(\theenumi)}

\item for each $\eta\in{\cal G}$ (the Lie algebra of $G$) and each
$h\in H$,
 $$\omega\bigl(\widehat \eta(h)\bigr)=\eta\,,$$
where $\widehat \eta$ is the fundamental vector field on $H$
associated to $\eta$;

\item for each $g\in G$, the pull-back $R_g^*\omega$ of the form
$\omega$ by the right translation $R_g:H\to H$, is
  $$R_g^*\omega=\Ad_{g^{-1}}\circ\, \omega\,.$$

\end{enum}

Conversely, any $1$-form on $H$ which takes its values in the Lie
algebra $\cal G$ of $G$, and which satisfies these two properties,
is the connection form of a $G_\gamma$-connection on the principal
bundle $H(B,G,G_\gamma,\overline H)$.
\end{proposition}

\section{Cartan connections seen by Ehresmann.} We will see in this section that connections on fibre bundles with a Lie group as
structure group offer a very nice setting for a rigorous
presentation of Cartan connections.

\begin{proposition}\label{propprincipalsubbundle}
Let $E(B,F,G,H)$ be a smooth bundle with a Lie group $G$ as structure group. We
assume that $F$ is an homogeneous space $G/G'$, with $G'$ a closed
subgroup of $G$. We denote by $o\in F$ the image of the unit element
$e\in G$ by the projection $\pi_G:G\to F=G/G'$. Moreover, we assume
that there exists a smooth section $s_0:B\to E$ of the bundle
$E(B,F,G,H)$. Let $H'$ be the subset of elements $h\in H$ such that
  $$h(o)=s_0\bigl(\pi_H(h)\bigr)\,.$$
Then $H'$ is the total space of a principal bundle
$H'(B,G',G'_\gamma,\overline{H'})$ with $G'$ as structure group.
\end{proposition}

\begin{dfn}\label{defcartanconnection}
Under the assumptions and with the notations of the previous
proposition, let $\omega$ be the connection form of a $G$-connection on $H$, and let $\omega_{H'}=i_{H'}^*\omega$ be the $1$-form
induced on $H'$ by $\omega$. We assume that for
any $h' \in H'$ and any vector $V'\in T_{h'}H'$, we have
 $$\omega_{H'}(V')=\omega(V')=0\quad\hbox{if and only if}\quad V'=0\,.$$
Moreover, we assume that
 $$\dim B=\dim F\,.$$
Such a connection on the bundle $E(B,G/G',G,H)$ (or on the
associated principal bundle $H(B,G,G_\gamma,\overline H)$) is called
a \emph{Cartan connection}.
\end{dfn}

\begin{remark}\label{reminducedformofcartanconnection}
The assumptions and notations are those of
Proposition~\ref{propprincipalsubbundle} and
Definition~\ref{defcartanconnection}.
Although defined on the total space $H'$ of the principal
bundle $H'(B,G',G'_\gamma,\overline{H'})$, the form $\omega_{H'}$ is
not a connection form on that bundle in the sense of Charles Ehresmann, 
since it takes its values in the larger Lie algebra $\cal G$, not in 
the Lie algebra ${\cal G}'$ of $G'$. 
However, the connection $1$-form $\omega$, defined on $H$, is
fully determined by the form $\omega_{H'}$ induced on $H'$: its
values on a vertical vector tangent to $H$ (maybe not tangent to
$H'$) at a point of $H'$ is known, since such a vector is the value
of a fundamental vector field associated to an element of the Lie
algebra $\cal G$; and its values at points of $H$ outside $H'$ can
be deduced from its values on $H'$ by using the formula
  $$R_g^*\omega=\Ad_{g^{-1}}\circ\,\omega\,.$$
Moreover, let $H'(B,G',G'_\gamma,\overline{H'})$ be a principal bundle
whose structure group $G'$ is a closed subgroup of a Lie group $G$, and let
$\omega_{H'}$ be a $1$-form defined on $H'$, taking its values in the Lie algebra
$\cal G$, satisfying the two properties:
  
\begin{enum}{(iii)}

\renewcommand\theenumi{\roman{enumi}}
\renewcommand\labelenumi{(\theenumi)}

\item for each $\eta'\in{\cal G}'$ (the Lie algebra of $G'$) and each
$h'\in H'$,
 $$\omega_{H'}\bigl(\widehat{\eta'}(h')\bigr)=\eta'\,,$$
where $\widehat {\eta'}$ is the fundamental vector field on $H'$
associated to $\eta'$;

\item for each $g'\in G'$, the pull-back $R_{g'}^*\omega_{H'}$ of the form
$\omega_{H'}$ by the right translation $R_{g'}:H'\to H'$, is
  $$R_{g'}^*\omega_{H'}=\Ad_{{g'}^{-1}}\circ\, \omega_{H'}\,.$$
  
\end{enum}  
  
Given such a $1$-form, called a \emph{Cartan connection form}, we can build a
principal bundle $H(B, G, G_\gamma, \overline H)$ with $G$ as structure group, and
a connection form $\omega$ on $H$, in such a way that $H'$ is a subset of $H$ and that
$\omega_{H'}$ is the form infuced on $H'$ by the connection form $\omega$.  We set $H=H'\times_{G'}G$, quotient of the product $H'\times G$ by the equivalence relation for which two elements $(h'_1,g_1)$ and $(h'_2,g_2)$ in $H'\times G$ are equivalent if there exists $g'\in G'$ such that  $h'_2=R_{g'}(h'_1)$ and $g_2=(g')^{-1}g$. We identify 
$H'$ with the set of equivalence classes of elements in $H'\times\{e\}$, and we define
$\omega$ on $H$ as indicated above.

\end{remark}

\begin{remark}\label{remgeneralizedcartan}
The first condition in the above definition,
  $$\omega(V')=0\quad\hbox{if and only if}\quad V'=0\,,$$
implies $\dim B\leq\dim F$, since $\dim H'=\dim B+\dim G'$ and $\dim
{\cal G}=\dim G=\dim F+\dim G'$. For a Cartan connection, the second
condition in the above definition imposes the equality, $\dim B=\dim
F$. When the first condition of the above definition is kept and the
second condition  is dropped, Ehresmann calls $\omega$ a \lq\lq
connexion de Cartan au sens large\rq\rq, {\it i.e.} a generalized
Cartan connection.
\end{remark}

\begin{remark}\label{remparallelizable}
The existence of a Cartan connection implies that $H'$ is
parallelizable, since $\omega_{H'}$ determines an isomorphism of
$TH'$ onto $H'\times{\cal G}$.
\end{remark}

The next proposition indicates another important consequence of the
existence of a Cartan connection: the tangent bundle $TB$ to the
base manifold is isomorphic to the bundle of vertical vectors,
tangent to the fibres of the bundle $E(B,F,G,H)$ along the image of
the global section $s_0$.

\begin{proposition}\label{propvectorbundleisomorphism}
The assumptions and notations are those of
Proposition~\ref{propprincipalsubbundle} and
Definition~\ref{defcartanconnection}. Let $x\in B$, $w\in T_xB$.
There exists a unique vector $\sigma(w)\in T_{s_0(x)}E_x$ such that,
for any $h'\in H'_x$ and $W'\in T_{h'}H'$ which satisfies
$T\pi_{H'}(W')=w$, we have
  $$\sigma(w)=T_oh'\circ T_e\pi_G\circ \omega(W')\,.$$
The map $\sigma:TB\to V_{s_0(B)}E$ so obtained is a vector bundle
isomorphism, called the soldering isomorphism, of the tangent bundle
$TB$ onto the vector bundle $V_{s_0(B)}E$ of vertical vectors,
tangent to the fibres of $E$ along the image $s_0(B)$ of the global
section $s_0$.
\end{proposition}

\Proof In the formula which defines $\sigma$ we have identified the
Lie algebra $\cal G$ with the space $T_eG$, tangent to the Lie group
$G$ at its unit element $e$. We have only to show that $T_oh'\circ
T_e\pi_G\circ \omega(W')$ depends only on $w=T\pi_{H'}(W')$, not on
the choices of $W'$ and $h'$. If we keep $h'$ unchanged and replace
$W'$ by $W'_1$, $W'_1-W'$ is the value at $h'$ of a fundamental
vector field associated to an element of the Lie algebra ${\cal
G}'$, therefore $\omega(W'_1)-\omega(W')\in{\cal G}'$ and
$T_e\pi_{G}\circ\omega(W'_1)=T_e\pi_{G}\circ\omega(W')$. Now if we
replace $h'$ by $h'_1=h'\circ g'$, with $g'\in G'$, we may replace
$W'$ by $W'_1=TR_{g'}(W')$, so
\begin{equation*}
\begin{split}
  T_oh'_1\circ T_e\pi_G\circ \omega(W'_1)
  &=T_o(h'\circ g')\circ T_e\pi_G\circ \omega\bigl(TR_{g'}(W')\bigr)\\
  &=T_oh'\circ TL_{g'}\circ T_e\pi_G\circ
  \Ad_{{g'}^{-1}}\circ\,\omega(W')\\
  &=T_oh'\circ T_e\pi_G\circ \omega(W')\,,
  \end{split}
  \end{equation*}
since $T_e\pi_G\circ\Ad_{{g'}^{-1}}\circ\,\omega(W')
=TL_{{g'}^{-1}}\circ\omega(W')$.
\endproof

\begin{remark}\label{remtangency}
For each point $x\in B$, the soldering isomorphism determined by a
Cartan connection on the bundle $E(B,G/G',G,H)$ allows us to
identify the space $T_xB$, tangent to the base $B$ at $x$, with the
space $T_{s_0(x)}E_x$, tangent to the fibre $E_x$ at $s_0(x)$. It is
in that precise sense that the soldering isomorphism allows us to
consider that for each $x\in B$, the base manifold $B$ is tangent to
the fibre $E_x$, the point $x\in B$ being in contact with the point
$s_0(x)\in E_x$.
\end{remark}

\begin{dfn}\label{defdevelopmentcartan}
The assumptions are those of Definition~\ref{defcartanconnection}.
Let $t\mapsto x(t)$ be a smooth path in the base manifold $B$. The
\emph{development} of that path in the fibre $E_{x(t_0)}$ is the
development (in the sense of Definition~\ref{defdevelopment}) of the
smooth path in $E$,
  $$t\mapsto s_0\bigl(x(t)\bigr)\,.$$
\end{dfn}

\begin{remark}\label{remdevelopmentcartan}
The development $t\mapsto y(t)$ of the path $t\mapsto x(t)$ in the
fibre $E_{x(t_0)}$ can be obtained as follows. Let $t\mapsto h'(t)$
be a smooth path in $H'$ such that, for each $t$,
 $$\pi_{H'}\bigl(h'(t)\bigr)=x(t)\,.$$
Such a path always exists (and is not unique). It is not horizontal
for the connection $\overline C$ (except when $t\mapsto x(t)$ is
constant) since $\omega_H$ does not vanish when applied to nonzero
vectors tangent to $H'$. Let $t\mapsto h(t)$ be the horizontal lift
in $H$ of the path $t\mapsto x(t)$, through the point
$h(t_0)=h'(t_0)$. There exists a unique path $t\mapsto
g(t)=h(t)^{-1}\circ h'(t)$ in $G$ such that $g(t_0)=e$ and that, for
each $t$,
  $$h'(t)=h(t)\circ g(t)\,.$$
An easy calculation shows that he development of $t\mapsto x(t)$ in
$E_{x(t_0)}$ is the path
  $$t\mapsto y(t)=h'(t_0)\circ g(t)(o)\,,$$
where $o\in G/G'$ is the image of the unit element $e$ by the
projection $\pi_G:G\to G/G'$.
\par\smallskip
Moreover, the map
  $$\frac{dx(t)}{dt}\Bigm|_{t=t_0}\mapsto
    \frac{dy(t)}{dt}\Bigm|_{t=t_0}$$
is the vector bundle isomorphism $\sigma$ of
Proposition~\ref{propvectorbundleisomorphism}. That property is in
agreement with the fact than when the fibre $E_{x(t_0)}$ is
considered as tangent to the base $B$ at $x(t_0)$, the tangent
vectors to the path $t\mapsto x(t)$ at $x(t_0)$ and to its
development at $s_0\bigl(x(t_0)\bigr) $ should be related by the
isomorphism $\sigma$.
\end{remark}

\section{Examples of Cartan connections}

\subsection{Homogeneous spaces.} Let $G'$ be a closed subgroup of a
Lie group $G$ such that the left action of $G$ on the homogeneous
space $F=G/G'$ is effective. We consider the trivial $G$-bundle
$E(B,F,G,H)$ with $B=F=G/G'$, $E=B\times F=(G/G')\times(G/G')$, the
fibres being the subsets $\{x\}\times F$ (in other words, the
projection $\pi_E:B\times F\to B$ is the first projection). The
associated principal $G$-bundle is the trivial bundle
$H(B,G,G_\gamma,\overline H)$ with $H=B\times G$, an element
$h=(x,g)\in H$ being considered as the diffeomorphism of the
standard fibre $F$ onto the fibre $E_x=\{x\}\times F$:
 $$y\mapsto \bigl(x,g(y)\bigr)\,.$$
The structure group $G$ acts on $H$ on the right by
  $$\bigl((x,g)\,,\,g'\bigr)\mapsto (x,gg')\,.$$
For the global section $s_0:B\to B\times F$, we choose the diagonal
map $x\mapsto (x,x)$. The principal $G'$-bundle
$H'(B,G',G'_\gamma,\overline{H'})$ has as its total space $H'$ the
subset of elements $h'=(x,g)\in H=B\times G$ which map the point
$o=\pi_G(e)$ of $F$ into the image of the global section $s_0$, {\it
i.e.} which are such that $x=g(o)$. Therefore $H'$ can be identified
with $G$ by means of $\bigl(g(o),g\bigr)\mapsto g$. The fibre of
$H'$ over a point $g(o)=\pi_G(g)$ of $B$ is the coset $gG'$.
\par\smallskip
On the trivial bundle $B\times F$, we consider the integrable
connection $C$ for which the parallel transport along any smooth
path in $B$ with end points $x_0$ and $x_1$ is the map
 $$\{x_0\}\times F\to\{x_1\}\times F\,,\quad (x_0,y)\mapsto
 (x_1,y)\,.$$
Given a smooth path $t\mapsto x(t)$ in $B$, its development in the
fibre $\bigl\{x(t_0)\}\times F$ is $t\mapsto \bigl(x(t_0),\,
x(t)\bigr)$. The connection considered here being integrable, the
soldering isomorphism integrates in the following sense: for each
point $x$ in the base manifold $B$, there is a natural
diffeomorphism $y\mapsto (x,y)$ of the base $B$ onto the fibre
$\{x\}\times F$, whose differential at $x$ is the soldering
isomorphism.
\par\smallskip
When we identify $H'$ with $G$ as indicated above, the $1$-form
$\omega_{H'}$ induced on $H'$ by the connection form is the left
Maurer-Cartan form on $G$ which, for any $g\in G$ and $X\in T_gG$,
yields $\omega_{H'}(X)=TL_{g^{-1}}(X)$; of course $\omega_{H'}(X)=0$
if and only if $X=0$, as required by
Definition~\ref{defcartanconnection}.

\subsection{Affine connections.} Let $M$ be a smooth $n$-dimensional
manifold. Its tangent bundle $TM$ is generally considered as a
vector bundle, with the linear group $\GL(n,\RR)$ as structure
group; but here we consider it as an \emph{affine} bundle, with the
group $\Aff(n,\RR)$ of affine transformations of $\RR^n$ as
structure group. We set $G=\Aff(n,\RR)$, $G'=\GL(n,\RR)$; the
standard fibre $G/G'$ is identified with $\RR^n$, considered as an
affine space. The corresponding principal bundle
$H(M,G,G_\gamma,\overline H)$ has as total space the set $H$ of
affine isomorphisms of $\RR^n$ onto the tangent space $T_xM$ at some
point $x\in M$. We take as global section $s_0$ of $TM$ the zero
section (for each $x\in M$, $s_0(x)$ is the zero vector at $x$). Let
$H'$ be the subset of elements in $H$ which map the origin of
$\RR^n$ into the image of $s_0$. In other words, an element $h'\in
H'$ is a \emph{linear} isomorphism of $\RR^n$ onto the tangent space
$T_xM$ (considered now as a vector space) at some point $x\in M$,
and $H'$ is the total space of the $G'$-principal bundle
$H'(M,G',G'_\gamma,\overline{H'})$.
\par\smallskip
Let $\overline C$ be a connection on the principal bundle
$H(M,G,G_\gamma,\overline H)$. Its connection form $\omega$ takes
its value in the Lie algebra of the affine group, which is a
semi-direct product ${\cal G}={\cal G'}\times \RR^n$ of the Lie
algebra $\cal G'$ of the linear group and of the Lie algebra $\RR^n$
of the group of translations in $\RR^n$. We may therefore write
 $$\omega=\omega'+\omega''\,,$$
where $\omega'$ takes its values in $\cal G'$ and $\omega''$ in
$\RR^n$. In the same way, the form $\omega_{H'}$ induced by $\omega$
on $H'$ can be written
  $$\omega_{H'}=\omega'_{H'}+\omega''_{H'}\,,$$
where $\omega'_{H'}$ takes its values in ${\cal G}'$ and
$\omega''_{H '}$ in $\RR^n$.
\par\smallskip
We see that $\omega'_{H'}$ is a connection form on the
$G'$-principal bundle $H'(M,G',G'_\gamma,\overline{H'})$; in other
words it is the connection form of a \emph{linear} connection.
\par\smallskip
The second term $\omega''_{H'}$ can be identified with a linear
endomorphism $\sigma:TM\to TM$ as follows. Let $x\in M$, $h'\in
H'_x$, $w\in T_xM$ and $W\in T_{h'} H'$ such that $T\pi_{H'}(W)=w$.
We set
  $$\sigma(w)=T_eh'\circ \omega''(W)\,.$$
As in the proof of Proposition~\ref{propprincipalsubbundle}, we see
that $\sigma(w)$ only depends on $w$, not on the choices of $h'$ and
$W$, and that $\sigma$ is a linear endomorphism of the tangent
bundle $TM$. Moreover, $\omega''_{H'}$ is fully determined by
$\sigma$.
\par\smallskip
Conversely, let $\varpi$ be the connection form of a linear
connection, {\it i.e.} a connection on the principal bundle
$H'(M,G',G'_\gamma,\overline H')$ of linear frames, and $\sigma$ be
a linear endomorphism of $TM$. There exists a unique connection on
the principal bundle $H(M,G,G_\gamma,\overline H)$ of affine frames
whose connection form $\omega=\omega'+\omega''$ is such that
$\omega'_{H'}=\varpi$ and that the the linear endomrphism of $TM$
determined by $\omega''_{H'}$ is $\sigma$. So a connection on the
principal bundle of affine frames can be identified with a pair made
by a linear connection and a linear endomorphism of $TM$.

\begin{remark} The Lie algebra of the affine group is a
semi-direct product ${\cal G}={\cal G'}\times\RR^n$ of the Lie
algebra $\cal G'$ of the linear group and the Lie algebra $\RR^n$ of
the group of translations $\RR^n$. Let $\pi_G:G\to G/G'$ be the
canonical projection. When $\cal G$ is identified with $T_eG$,
$G/G'$ and $T_o(G/G')$ with $\RR^n$, the second projection ${\cal
G}\to\RR^n$ is simply the tangent map $T_e\pi_G$ at the unit
element.
\end{remark}

The connection form $\omega$ on the frame of affine bundles is a
Cartan connection in the sense of
Definition~\ref{defcartanconnection} if and only if $\omega_{H'}$
has zero kernel, {\it i.e.} for each $h'\in H'$, $\overline
C_{h'}\cap T_{h'}H'=\{0\}$. Let $S'_{h'}$ and $S''_{h'}$ be the
subspaces of $T_{h'}H$ spanned by the values at $h'$ of fundamental
vector fields associated, respectively, to elements in the Lie
subalgebras $\cal G'$ and $\RR^n$ of $\cal G$. We have the direct
sum decompositions
  $$\ker T_{h'}\pi_H=S'_{h'}\oplus S''_{h'}\,,\quad
  T_{h'}H=S'_{h'}\oplus S''_{h'}\oplus \overline C_{h'}\,.$$
Since $\ker T_{h'}\pi_{H'}=T_{h'}H'\cap \ker T_{h'}\pi_H=S'_{h'}$ we
see that $\overline C_{h'}\cap T_{h'}H'=\{0\}$ if and only if
$T_{h'}H'\cap (S''_{h'}\oplus \overline C_{h'})$ is the graph of a
linear isomorphism of $C_{h'}$ onto $S''_{h'}$. Using the above
remark, we conclude that the connection form $\omega$ is a Cartan
connection if and only if the endomorphism $\sigma:TM\to TM$
determined by $\omega$ is an isomorphism. When that condition is
fulfilled, $\sigma$ is the soldering isomorphism. It is very natural
to impose to $\sigma$ to be the identity isomorphism of $TM$, since
that ensures that any smooth path $t\mapsto x(t)$ in $M$ has, at a
point $x(t_0)$, the same tangent vector as its development in
$T_{x(t_0)}M$. Affine connections considered by
\'E.~Cartan~\cite{cartan1} satisfy that property. Contrary to
Ehresmann~\cite{ehresmann2} and Lichnerowicz~\cite{lichnerowicz} who
call \emph{affine connection} any connection on the principal bundle
of affine frames, Kobayashi and Nomizu~\cite{kobayashi1} reserve
that name for a Cartan connection whose soldering isomorphism is the
identity of $TM$, and speak of \emph{generalized affine connections}
for all other connections on the principal bundle of affine frames.

\subsection{Projective connections.}
In this subsection the standard fibre $F=\P(n,\RR)$ is the standard
$n$-dimensional real projective space, quotient of $\RR^{n+1}-\{0\}$
by the equivalence relation of colinearity. The structure group will
be the corresponding projective linear group $G=\PGL(n,\RR)$,
quotient of $\GL(n+1,\RR)$ by the group of homotheties. We denote by
$x_0,x_1,\ldots,x_n$ the standard coordinates in $\RR^{n+1}$. Let
$D$ be the straight line in $\RR^{n+1}$ defined by the equations
$x_1=\cdots=x_n=0$, $G'$ be the closed subgroup of $G$ made by
equivalence classes of elements $g\in \GL(n+1,\RR)$ which map $D$
onto itself, and $o\in F$ be the image of $D-\{0\}$ by the
projection $\RR^{n+1}-\{0\}\to \P(n,\RR)$. The standard fibre
$F=\P(n,\RR)$ will be identified with the homogeneous space $G/G'$.
For each $g\in \GL(n,\RR)$, let $\widehat g\in\PGL(n,\RR)$ be the
map
 $$\widehat g:\P(n,\RR)\to\P(n,\RR)\,,\quad
  [x_0,x_1,\ldots,x_n]\mapsto[x_0,g(x_1,\ldots,x_n)]\,.$$
The injective map $g\mapsto\widehat g$ will allow us to identify
$\GL(n,\RR)$ with a closed subgroup $G''$ of $\PGL(n,\RR)$,
contained in the subgroup $G'$ defined above.We observe that $G'$ is
the subgroup made by elements in $G$ which leave fixed the point
$o\in F$, and that $G''$ is the subgroup made by elements which, in
addition, leave globally invariant the hyperplane at infinity (image
of the hyperplane defined by the equation $x_0=0$ by the projection
$\bigl(\RR^{n+1}-\{0\}\bigr)\to \P(n,\RR)$).

\par\smallskip Now let $M$ be an
$n$-dimensional smooth manifold. We consider the direct sum of the
trivial bundle $M\times \RR$ and of the tangent bundle $TM$. We take
the complementary part of the image of the zero section in the total
space of that bundle and quotient it by the equivalence relation of
colinearity in the fibres. We obtain a fibre bundle $E(M,F)$ whose
standard fibre is the projective space $F=\P(n,\RR)$. The tangent
bundle can be considered as a dense subset of the total space $E$ of
that bundle, if we identify, for each $x\in M$ and $v\in T_xM$, the
vector $v$ with the element in $E_x$, equivalence class of $(1,v)$.
Since the tangent bundle $TM$ and the trivial bundle $M\times\RR$
admit as structure group, respectively, the linear group $GL(n,\RR)$
and the trivial group $\{e\}$, and since we have identified
$\GL(n,\RR)$ with a closed subgroup of $\PGL(n,\RR)$, the bundle
$E(M,F)$ admits $G=\PGL(n,\RR)$ as structure group. Following
Ehresmann's notations, we will denote it by $E(M,F,G,H)$, where $H$
is the principal bundle of projective frames in the fibres of
$E(M,F)$.
\par\smallskip
For each $x\in M$, we denote by $s_0(x)$ the equivalence class of
$(1,0_x)$, where $0_x$ is the zero vector at $x$. We may now define
the subset $H'$ of elements $h'\in H$ which map the element $o\in F$
into the image of $s_0$, and observe that it is the total space of a
principal bundle $H'(M,G',G'_\gamma,\overline{H'})$. A Cartan
connection, in the sense of Definition~\ref{defcartanconnection}, is
a connection form $\omega$ on the principal bundle
$H(M,G,G_\gamma,\overline H)$ such that the form $\omega_{H'}$
induced on $H'$ has zero kernel. The soldering isomorphism can still
be considered as an isomorphism of the tangent bundle $TM$, since we
have identified $TM$ with an open dense subset of $E$. A
\emph{projective connection} is a Cartan connection whose soldering
isomorphism is the identity of $TM$. Cartan and Ehresmann have shown
that projective connections exist on any manifold $M$.

\subsection{Conformal connections.}

We follow the presentation of Kobayashi~\cite{kobayashi2}. On  the
space $\RR^{n+2}$, with coordinates $(x_0,x_1,\ldots,x_{n+1})$, let
$Q$ be the quadratic form
  $$Q(x_0,\ldots,x_{n+1})=x_1^2+\cdots+x_{n+1}^2-x_0^2\,.
  $$
Its signature is $(n+1,1)$. Let $\O(n+1,1)$ be the subgroup of
$\GL(n+2,\RR)$ of elements which leave the quadratic form $Q$
unchanged, and let $Q_0$ be the cone, subset of $\RR^{n+2}$, defined
by the equation $Q=0$. Let us call \emph{rays} the straight lines
through the origin in $\RR^{n+2}$. Any element $g\in\O(n+1,1)$ maps
a ray onto another ray, and a ray contained on $Q_0$ on another ray
contained in $Q_0$. The set of rays is the projective space
$\P(n+1,\RR)$. Therefore the group $O(n+1,1)$ acts on $\P(n+1,\RR)$
and that action leaves invariant the image $\cal M$ of $Q_0$ by the
projection $\bigl(\RR^{n+2}-\{0\}\bigr)\to\P(n+1,\RR)$. That action
is transitive on $\cal M$. We denote by $o\in{\cal M}$ the ray
defined by the equations $x_0-x_{n+1}=0\,,\ x_1=\cdots=x_n=0$, and by
$G'\subset\O(n+1,1)$ its stabilizer. We will identify $\cal M$,
which is called the \emph{M\"{o}bius space}, with $\O(n+1,1)/G'$.
\par\smallskip

Let $\Pi$ be the affine hyperplane in $\RR^{n+2}$ defined by the
equation $x_0=1$. The affine injective map, defined on
$\RR^{n+1}$, with values in $\RR^{n+2}$,
  $$(y_1,\ldots,y_{n+1})\mapsto (x_0,x_1,\ldots,x_{n+1})\,,
  $$
with
  $$x_0=1\,,\ x_1=y_1\,,\ \ldots,\ x_n=y_n\,,\
  x_{n+1}=y_{n+1}\,,
  $$
has $\Pi$ as image and maps the sphere
  $$S^n=\left\{\,(y_1,\ldots,y_{n+1})\in\RR^{n+1};
  \sum_{i=1}^{n+1}y_i^2=1\,\right\}
  $$
onto $Q_0\cap \Pi$. Since $\Pi$ meets each ray contained in $Q_0$ at
a unique point, $Q_0\cap \Pi$ is diffeomorphic to $\cal M$.
Therefore, by composition with the projection
$\bigl(\RR^{n+2}-\{0\}\bigr)\to\P(n+1,\RR)$, we obtain a
diffeomorphism of $S^n$ onto the M\"{o}bius space $\cal M$.
\par\smallskip

Let $\Pi_1$ be the affine hyperplane in $\RR^{n+2}$ defined by the
equation $x_0+x_{n+1}=1$. It meets each ray contained in $Q_0$ at a unique
point, except the ray defined by the equations
$x_0+x_{n+1}=0\,,\ x_1=\cdots=x_n=0$. The smooth injective map, defined on $\RR^n$
and with values in $\RR^{n+2}$,
  $$(z_1,\ldots,z_{n})\mapsto (x_0,x_1,\ldots,x_{n+1})\,,
  $$
with
  $$x_0=\frac{1+\sum_{i=1}^nz_i^2}{2}\,,\ 
  x_1=z_1\,,\ \ldots,\ x_n=z_n\,,\
  x_{n+1}=\frac{1-\sum_{i=1}^nz_i^2}{2}\,,
  $$
has $Q_0\cap \Pi_1$ as image. Composed with the projection
$\bigl(\RR^{n+2}-\{0\}\bigr)\to\P(n+1,\RR)$, that map yields a
smooth injective map of $\RR^n$ into the M\"{o}bius space $\cal M$, whose
image is ${\cal M}$ minus one point, the missing point corresponding
to the ray defined by the equations $x_0+x_{n+1}=0\,,\ x_1=\cdots=x_n=0$. The
image, by that map, of the origin of $\RR^n$ is the point $o\in{\cal
M}$. Moreover, that map is equivariant with respect to the actions
of the orthogonal group $\O(n)$, on $\RR^n$ and on $\RR^{n+2}$, with
the convention that an element $g\in \O(n)$ is identified with the
element $\widehat g\in \O(n+1,1)$ which maps
$(x_0,x_1,\ldots,x_n,x_{n+1})$ on
$\bigl(x_0,g(x_1,\ldots,x_n),x_{n+1}\bigr)$.

\begin{remark}
When the M\"{o}bius space $\cal M$ is identified with the sphere $S^n$,
the injective map $\RR^n\to {\cal M}$ defined above is the inverse
of the stereographic projection from $S^n$ minus its south pole onto $\RR^n$.
\end{remark}

\par\smallskip

Now let $(M,g)$ be an $n$-dimensional smooth Riemannian manifold. We
consider the direct sum of two copies of the trivial bundle $M\times
\RR$ and of the tangent bundle $TM$, in the following order:
$(M\times\RR)\oplus TM\oplus (M\times\RR)$. To shorten the notation,
we will denote that bundle by $\RR\oplus TM\oplus \RR$. Each element
of the total space of that bundle is a triple $(v_0,v,v_{n+1})$,
with $v_0$ and $v_{n+1}\in \RR$ and $v\in TM$, its projection on the
base $M$ being the projection of $v$. Let $Q$ be the quadratic form,
defined on the fibres of that bundle,
  $$Q(v_0,v,v_{n+1})=g(v,v)+v_{n+1}^2-v_0^2\,,$$
and $Q_0$ be the subspace of $\RR\oplus TM\oplus\RR$ defined by the
equation $Q=0$. We consider the projective bundle $\P(\RR\oplus
TM\oplus\RR)$, quotient of $\bigl(\RR\oplus TM\oplus\RR-\{0\}\bigr)$
by the equivalence relation of colinearity in the fibres. We have
denoted by $\{0\}$ the image of the zero section in the bundle
$\RR\oplus TM\oplus\RR$. The image of $(Q_0-\{0\})$ by the
projection of $\bigl(\RR\oplus TM\oplus\RR-\{0\}\bigr)$ onto
$\P(\RR\oplus TM\oplus\RR)$ is the total space of a fibre bundle
$E(M,{\cal M})$, with base the manifold $M$, and with standard fibre
the M\"{o}bius space $\cal M$. The tangent bundle $TM$ can be considered
as an open, dense subset of $E$ by identifying each vector $v\in TM$
with the equivalence class of 
 $\Bigl(\bigl(1+g(v,v)\bigr)/2,\,v,\,\bigl(1-g(v,v)\bigr)/2\Bigr)$. 
Since
the manifold $M$ is equipped with a Riemannian metric, the tangent
bundle $TM$ admits $\O(n)$ as structure group. This is true for the
bundle $E(M,{\cal M})$ too. But as seen above, $\O(n)$ can be 
considered as a subgroup of $\O(n+1,1)$, which acts on the standard
fibre $\cal M$. The bundle $E(M,{\cal M})$ therefore has
$G=\O(n+1,1)$ as structure group and, using Ehresmann's notation, we
will denote it by $E(M,{\cal M},G,H)$. The total space $H$ of the
corresponding principal bundle is the space of \emph{conformal
frames} on the Riemannian manifold $(M,g)$.

\par\smallskip
For each $x\in M$, we denote by $s_0(x)$ the equivalence class of
$(1,0_x,1)$, where $0_x$ is the zero vector at $x$. The subset $H'$
of elements $h'\in H$ which map the element $o\in{\cal M}$ into the
image of $s_0$ is the total space of a principal bundle
$H'(M,G',G'_\gamma,\overline{H'})$. A Cartan connection, in the
sense of Definition~\ref{defcartanconnection}, is a connection form
$\omega$ on the principal bundle $H(M,G,G_\gamma,\overline H)$ such
that the form $\omega_{H'}$ induced on $H'$ has zero kernel. The
soldering isomorphism can still be considered as an isomorphism of
the tangent bundle $TM$, since we have identified $TM$ with an open
dense subset of $E$. A \emph{conformal connection} is a Cartan
connection whose soldering isomorphism is the identity of $TM$.
Cartan and Ehresmann have shown that conformal connections exist on
any manifold $M$.

\section{Applications of connections.}

\subsection{Gravitation.}
The first very important application of the notion of connection is
probably to be found in the theory of General Relativity, in which
the Levi-Civita connection associated to the pseudo-Riemannian
structure on space-time plays a key role. According to the title of
his paper~\cite{cartan1}, \lq\lq Sur les espaces \`{a} connexion affine
et la th\'{e}orie de la relativit\'{e} g\'{e}n\'{e}ralis\'{e}e\rq\rq, \'E.~Cartan was
probably, for a large part, motivated by possible physical
applications when he investigated the properties of connections. He
explains how, in the framework of classical, non relativistic
mechanics, a gravitation field ({\it i.e.} an acceleration field)
can be included in the geometry of space-time by the use of an
appropriate affine connection. We present here this idea in the more
modern language of Ehresmann. Our manifold $M$ is the non
relativistic space-time. For simplicity we assume that $\dim M=2$
(we take into account only one dimension for space). The choice of a
Galilean frame and of units for time and length allows us to
identify $M$ with $\RR^2$, with coordinates $(t,x)$. The Lie group
$G$, which will be called the \emph{affine Galileo group}, is the
group of affine transformations of $\RR^2$ of the form
  $$g:\RR^2\to \RR^2\,,\quad (t,\,x)\mapsto (t'=t+a,\, x'=x+b+vt)\,.$$
An element $g$ of $G$ is therefore a triple $(v,\,a,\,b)$, and the
composition law of $G$ is
  $$(v_2,\,a_2,\,b_2)(v_1,\,a_1,\,b_1)=(v_2+v_1,\, a_2+a_1,\,
  b_2+b_1+v_2a_1)\,.
  $$
Let $G'$ be the linear Galileo group, {\it i.e.} the subgroup of
elements of $G$ of the form $(v,0,0)$. The homogeneous space
$F=G/G'$ can be identified with $\RR^2$, with coordinates $(a,b)$.
The bundle $E(M,F,G,H)$ is simply the tangent bundle $TM$,
considered as an affine bundle, its structure group being restricted
to the affine Galileo group (instead of the full affine group). The
corresponding principal bundle $H(M,G,G_\gamma,\overline H)$ is the
bundle of affine Galilean frames in $TM$. Its total space is simply
the product $H=M\times G$. The subset of elements in $H$ which map
the origin of $\RR^2$ in the image of the zero section of $TM$ is
$H'=M\times G'$. It is the total space of a principal bundle
$H'(M,G',G'_\gamma,\overline{H'})$.
\par\smallskip
An affine Galilean connection is determined by a connection form
$\omega$ on $H$, with values in the Lie algebra $\cal G$ of $G$,
which induces on $H'$ a form $\omega_{H'}$ with zero kernel. In
addition we impose to that form to be such that the corresponding
soldering isomorphism is the identity of $TM$. The Lie algebra $\cal
G$ can be identified with $\RR^3$, by means of the basis
$(\varepsilon_v,\varepsilon_a,\varepsilon_b)$ which corresponds to
the coordinates $(v,a,b)$ on $G$. Taking into account the
equivariance properties of $\omega$ and the fact that the
corresponding soldering isomorphism is the identity of $TM$, we
obtain
  $$\omega=(-V\,dt-W\,dx+dv)\varepsilon_v +(dt+da)\varepsilon_a
  +\bigl(-(v+aV)\,dt+(1-aW)\,dx-v\,da+db\bigr)\varepsilon_b\,,
  $$
where $(t,x,v,a,b)$ are the coordinates on $H=M\times G$, $V$ and
$W$ being two smooth functions on $M$ (therefore depending only on
the coordinates $(t,x)$).

\par\smallskip

The development of the trajectory $t\mapsto \bigl(t,\,x(t)\bigr)$ of
a particle is a straight line if and only if the functions $V$ and
$W$, on which the connection form depend, satisfy
  $$V\bigl(t,x(t)\bigr)+W\bigl(t,x(t)\bigr)
  \frac{dx(t)}{dt}-\frac{d^2x(t)}{dt^2}=0\,.
  $$
We assume now that there exists on $M$ a gravity field $g$ (which
may eventually depend on time $t$ and space location $x$). The
equation of motion of a material particle submitted to that
gravitational field is
  $$\frac{d^2x(t)}{dt^2}=g(t,x)\,.
  $$
Therefore, if we choose $W=0$ and $V(t,x)=g(t,x)$, the development
of the trajectory of any material particle submitted to the gravity
field $g$, but to no other forces, is a straight line. When $g$ is a
constant, the corresponding connection is integrable: a nonlinear
change of coordinates in space-time eliminates the gravity force (in
agreement with Einstein's thought experiment in which an observer in
a lift in free fall no more feels the gravity force). This is no
more true when $g$ is not constant. It would be interesting to
develop the example in which $g$ is the acceleration field in
Kepler's problem (in a non relativistic space-time of dimension 3);
for a suitably chosen connection, the developments of the Keplerian
trajectories of the planets should be straight lines.

\subsection{Geometric quantization.}
Let $(M,\Omega)$ be a symplectic manifold. A \emph{prequantization}
of that symplectic manifold is a principal bundle $P(M,S^1,S^1,H)$
with base $M$ and with the circle $S^1$ as structure group, endowed
with a connection $1$-form $\omega$ whose curvature is $\Omega$.
\par\smallskip

Since the Lie algebra of $S^1$ can be identified with $\RR$, the
connection form $\omega$ can be considered as a \emph{contact form}
on $P$. According to a theorem by B.~Kostant~\cite{kostant} and
J.-M.~Souriau~\cite{souriau}, there exists a prequantization of
$(M,\Omega)$ if and only if the cohomology class of $\Omega$ is
integer.

\subsection{Phases in Mechanics.}

Various uses of connections are made in the mathematical description
of mechanical systems.
\par\smallskip

Let us consider a Hamiltonian system, depending of some parameters,
which for any fixed value of these parameters, is completely
integrable. The motion of the system, for a fixed value of the
parameter, is quasi-periodic on a Lagrangian torus of phase space.
At a certain time, the parameters vary slowly, describe a closed
loop in the space of values of the parameters, and after taking
again their initial values, remain constant. The motion of the
system becomes again quasi-periodic on the same Lagrangian torus,
but with a change of phase (the Hannay and Berry phase). This change
of phase is interpreted as the holonomy of an Ehresmann connection
in the works of Marsden, Montgomery and Ratiu~\cite{marsden,
montgomery}.

\subsection{Nonholonomic constraints.}
 Several different approaches have been used for
 the mathematical description of mechanical systems with
 constraints. In one of these approaches, the configuration space
 of the system is a smooth manifold and the constraints are
 described by a vector (or sometimes an affine) sub-bundle $C$ of the
 tangent bundle $TQ$. The admissible motions of the system are
 smooth curves $t\mapsto x(t)$ in $Q$ such that, for all $t$,
  $$\frac{dx(t)}{dt}\in C_{x(t)}\,.$$
\par\smallskip

J.~Koiller~\cite{koiller} considered systems where the configuration
space $Q$ is the total space of a principal bundle over a base $B$,
with a Lie group $G$ as structure group, and where the constraint
$C$ is a connection on that principal bundle.
\par\smallskip

\subsection{Active constraints.}
Let us consider a mechanical system in which some geometric
constraints can be acted on, as a function of time, in order to
control the motion of the system. For example, a cat in free fall
can change the shape of her body to try (generally with success) to
reach the ground on her feet.
\par\smallskip
For the mathematical description of such systems~\cite{marle}, we
use a manifold $Q$ as configuration space, and a smooth submersion
$\pi:Q\to S$ onto another manifold $S$ (the space of shapes of the
cat's body, or more generally the space of possible states of the
active constraint).
\par\smallskip
The dynamical properties of the system (other than those used to
change the value of the active constraint) are described by a
Lagrangian $L:TQ\to\RR$. Let $V^*Q$ be the dual bundle of the
vertical sub-bundle $VQ=\ker T\pi\subset TQ$. We can identify
$V^*Q$ with the quotient bundle $T^*Q/(VQ)^0$. Let
$\zeta:T^*Q\to V^*Q$, $q:V^*Q\to Q$ and $\widetilde \pi=\pi\circ
q:V^*Q\to S$ be the projections. When the Lagrangian $L$ is
 $$L(v)=\frac{1}{2}g(v,v)-P(x)\,,\quad\hbox{with}\quad x\in Q\,,\
 v\in T_xQ\,,$$
there is on the bundle $V^*Q\to S$ an Ehresmann connection (called
the dynamical connection) which can be used to determine the way in
which an infinitesimal change of the state of the active constraint,
represented by a vector tangent to $S$, affects the motion of the
mechanical system.

\par\smallskip
Let us call \emph{kinetic connection} the Ehresmann connection, on
the bundle $\pi:Q\to S$, for which the horizontal lift at $x\in Q$
of a vector $v\in T_{\pi(x)}S$ is the unique $w\in T_xQ$, orthogonal
(with respect to $g$) to the vertical subspace $\ker T_x\pi$,  such
that $T\pi(w)=v$. The dynamic connection is characterized by the two
properties:

\begin{enum}{(iii)}
\renewcommand\theenumi{\roman{enumi}}
\renewcommand\labelenumi{(\theenumi)}

\item the horizontal lift at $z\in \widetilde\pi^{-1}(s)$ of a vector
$v\in T_sS$ with respect to the dynamical connection projects on $Q$
onto the horizontal lift at $x=q(z)$ of $v$ with respect to the
kinetic connection;

\item the horizontal lift to $V^*Q$ of any smooth vector field on
$S$, with respect to the dynamical connection, is an infinitesimal
automorphism of the Poisson structure of $V^*Q$.
\end{enum}
\par\smallskip
\subsection{Maxwell's equations.}
The famous Maxwell's equations are usually written (see for
example~\cite{guillemin})
 \begin{equation*}
\begin{aligned}
\rot\overrightarrow E+\frac{\partial\overrightarrow B}{\partial
t}&=0\,, &\rot\overrightarrow H-\frac{\partial\overrightarrow
D}{\partial t}&=4\pi \overrightarrow j\,,\\
\div\overrightarrow B&=0\,, &\div \overrightarrow D&=4\pi\rho\,,
\end{aligned}
\end{equation*}
where $\overrightarrow E$ is the electric field, $\overrightarrow B$
the magnetic induction, $\overrightarrow D$ the diplacement current,
$\overrightarrow H$ the magnetic field, $\rho$ the electric charge
density and $\overrightarrow j$ the current density. Moreover there
are constitutive equations which link $\overrightarrow E$ and
$\overrightarrow D$, $\overrightarrow B$ and $\overrightarrow H$,
 $$\overrightarrow D=\varepsilon_0\overrightarrow E\,,\quad
   \overrightarrow H=\frac{1}{\mu_0}\overrightarrow B\,.$$
\par\smallskip
Let us introduce the $2$-form on space-time (in which
$(t,\,x^1,\,x^2,\,x^3)$ are the coordinates, with respect to some
Galilean reference frame, some chosen units of time and length and
some orthonormal frame in space)
\begin{equation*}
  F
  =B_1\, dx^2\wedge dx^3+B_2\,dx^3\wedge dx^1+B_3\,dx^1\wedge dx^2 
  +(E_1\,dx^1+E_2\,dx^2+E_3\,dx^3)\wedge dt\,,
\end{equation*}
and
\begin{equation*}
  G
  =D_1\, dx^2\wedge dx^3+D_2\,dx^3\wedge dx^1+D_3\,dx^1\wedge dx^2
  -(H_1\,dx^1+H_2\,dx^2+H_3\,dx^3)\wedge dt\,.
\end{equation*}
Let us set
  $$J=\rho\, dx^1\wedge dx^2\wedge
  dx^3-(j_1\,dx^2\wedge dx^3+j_2\,dx^3\wedge dx^1+j_3\,dx^1\wedge dx^2)
  \wedge dt\,.$$

\par\smallskip
Then we have
 $$G=\sqrt{\frac{\varepsilon_0}{\mu_0}} * F\,,$$
where $*$ is the Hodge operator on the $4$-dimensional
pseudo-Riemannian manifold space-time.
 \par\smallskip
The velocity of light in vacuum is
 $$c=\sqrt{\frac{1}{\varepsilon_0\mu_0}}\,.$$
 \par\smallskip
Maxwell's equations become
  $$dF=0\,,\quad dG=4\pi J\,,\quad\hbox{with}\quad
  G=\sqrt{\frac{\varepsilon_0}{\mu_0}} * F\,.$$

\par\smallskip
Maxwell's equations can be made even more beautiful: the
electromagnetic $2$-form $F$ on space-time can be considered as the
curvature form of a connection on a principal bundle, with
space-time as base and the circle $S^1=\U(1)$ as structure group
(this idea was first introduced by H.~Weyl~\cite{weyl1}). The
$2$-form $F$ should then be considered as taking its values in the
Lie algebra $u(1)$ of $S_1$. The connection form $A$ such that
$F=DA$, is not unique: we may add a closed $1$-form (gauge
transformation). The first Maxwell's equation, $dF=0$, is
automatically satisfied. The second Maxwell's equation becomes
  $$D(*DA)=4\pi\sqrt{\frac{\mu_0}{\varepsilon_0}}J\,,$$
where $D$ is the covariant exterior differential operator.

\subsection{ Yang-Mills fields.}
Gauge theories generalize Maxwell's theory of electromagnetism
written in terms of connections. They use a principal bundle with
space-time as base and a non-Abelian group as structure group
($\U(1)\times \SU(2)\times \SU(3)$ in the so-called \emph{standard
model}).
\par\smallskip
They introduce a connection $1$-form $A$ on that bundle and lead,
for the curvature form $F$ of the connection, to field equations
similar to Maxwell's equations:
   $$ D_AA=F\,,
   \quad D_AF=0\quad\hbox{ (Bianchi identity)}\,,
   $$
 and
  $$D_A * F=J\,,$$
where $D_A$ is the covariant differential with respect to the
connection $A$ and $J$ a \lq\lq current\rq\rq\ which generalizes the
$4$-dimensional current density of Maxwell's theory.
\par\smallskip

\section{Acknowledgements} Many thanks to Paulette Libermann, for her very
careful reading of this text and her many suggestions for improvements. And again my warmest thanks to the organizers and supporting institutions of this great International Conference.

\end{document}